\documentclass[12pt]{article}

\newtheorem{theo}{Theorem}[section]
\newtheorem{lem}[theo]{Lemma}

\usepackage{graphics}
\input{epsf}

\def\C{{\rm\kern.24em \vrule width.02em height1.4ex
depth-.05ex\kern-.26em C}}

\def\be{\begin{equation}}
\def\ee{\end{equation}}
\def\beq{\begin{equation}}
\def\eeq{\end{equation}}

\usepackage{graphics}
\input{epsf}
 \date{}
\begin{document}
  \title{
  Zero-Energy Self-Similar  Solutions Describing  Singularity Formation In The   Nonlinear Schrodinger Equation   
    In Dimension $N=3$ }
  \author{William  C. Troy
     }
 
\maketitle

 \begin{abstract} 
  In dimension $N=3$  the cubic nonlinear Schrodinger (NLS)  equation has solutions which become singular, i.e. 
at a spatial point they blow up to infinity   in finite time.
  In 1972 Zakharov famously investigated finite time singularity formation in   
 the cubic nonlinear Schrodinger equation   as a model for
  spatial   collapse 
of   Langmuir   waves in plasma, the most abundant form of observed  matter in the universe. 
  Zakharov assumed that    (NLS)     blow up of solutions  
  is self-similar and radially symmetric, and  that 
singularity formation  
   can be modeled    by a solution of an associated   self-similar, complex  ordinary differential equation~(ODE).
  A   parameter $a>0$ appears in the ODE, and the dependent variable, 
   $Q,$   satisfies   $(Q(0),Q'(0))=(Q_{0},0),$ where $Q_{0}>0.$ 
    A fundamentally important step towards putting the Zakharov  model on a firm mathematical footing     is to  
   prove, when $N=3,$ 
     whether   values $a>0$ and $Q_{0}>0$ exist such that   $Q$ also satisfies the   
    physically important `zero-energy' integral constraint.   
   Since 1972 this  has remained an open problem.  Here, we resolve this issue  by proving  that   
    for {\sl every}   $a>0$  and $Q_{0}>0,$   
   $Q$      satisfies the 
    the  `zero-energy' integral constraint. 
    \end{abstract}




\section{Introduction } \label{sec:intro}
\setcounter{equation}{0} 
  The nonlinear Schrodinger
 system   
    \begin{eqnarray}
{\rm i}\psi_{t}+\Delta \psi+\psi|\psi|^{2}=0, ~~t>0, 
 \label{eq:odre55}
 \end{eqnarray}
     \begin{eqnarray}
\psi(x,0)=u_{0}(x),~~x \in R^{N} 
 \label{eq:odre55a}
 \end{eqnarray}
      is  often considered to be the
     simplest model of singularity formation in nonlinear dispersive systems~\cite{kop,new}. 
When dimension  $N \in  [2,4)$ there exists a wide class of initial conditions for which the solution of (\ref{eq:odre55})-(\ref{eq:odre55a}) 
 forms a singularity  at 
   time $T>0,$ i.e. as  $ t \to T^{-}$  the solution
 becomes infinite at a single spatial  point where   a growing and increasingly narrow peak forms~\cite{kop, ras}. 
When $N=2 $ problem (\ref{eq:odre55})-(\ref{eq:odre55a}) asises in modeling singularity formation in nonlinear optical media,   
     and finite time blowup of a solution   corresponds to an   extreme increase in field amplitude due to `self-focusing'~\cite{chai, ha, kop}.
   In dimension $N=3$ problem (\ref{eq:odre55})-(\ref{eq:odre55a}) 
    arises as the subsonic limit of the    Zakharov model for   Langmuir waves in plasma,  and 
     singularity formation  is  referred to as wave  `collapse'~\cite{kop,zh1972,z}.
 Plasma,  often described as the fourth state of matter, 
      consists mainly  of  charged particles~(ions) and/or electrons, and is   the   most abundant form  of observed matter 
     in the universe~\cite{chu}. On August 22, 1879 the existence of plasma was first reported   
      by Sir William Crookes,
      who identified it  
      as ``radiant matter''
      in a  lecture to the British Association for the Advancement of Science. In 1928  Langmuir  introduced  the word plasma  in his studies of
     plasma waves, i.e.    oscillations in the  density of `ionized gas'~\cite{langm}. 
  Over the general range $ 2 < N <4,$ multiple  investigations of (\ref{eq:odre55})-(\ref{eq:odre55a}) have led to new understandings of
       singularity formation both from 
  the numerical and analytical point 
  of view~\cite{adk, budd, fib96, fp98, kop, Land,  L1, L, MCL, new, rott, rott2, sul, ts92,z}. 
   In particular, when $N=3,$ the 1972 and 1984 investigations by  Zakharhov~\cite{zh1972,z}, the 1986 numerical study by McGlaughlin
    et al~\cite{MCL}, and the 
   1988 numerical  investigations of LeMesurier et al~\cite{L1, L}  and Landman et al~\cite{Land},   
 demonstrated that, for both symmetric and asymmetric initial data,  singularity formation occurs in
  a  spherically symmetric~(i.e. $r=|x|$)  and self-similar
  manner. Thus, for   $ N \in (2,4),$  these authors   
modeled singularity formation  at $r=0$ and time $t=T,$ 
by an exact   
solution   of the form 
    \begin{eqnarray}
 \psi = { \frac  {1} { {\sqrt { 2a(T-t)}}}} 
 \exp \left (    i\tilde \theta +
 {\frac {i}{2a}}    \ln \left (  { \frac    { T }   { T-t }    }                         \right ) \right )
 Q \left ( {\frac {r}{  {\sqrt {  2a(T-t)}} } }\right ), 
 \label{eq:odre55aa}
 \end{eqnarray}
  where $a>0$  and $\tilde \theta $~(fixed phase shift) are constants.    
 Let 
  $\zeta= {\frac {r}{  {\sqrt {  2a(T-t)}} } }.$ 
  Then  $\zeta \to \infty$ as $t \to T$
 from below, and $Q(\zeta)$ 
 solves  the profile equation 
  \begin{eqnarray}
 Q''+{\frac {N-1}{\zeta}}Q'-Q+ia(Q +\zeta Q')+Q|Q|^{2}=0,~~\zeta \ge 0, 
   \label{eq:odreaa}
 \end{eqnarray}
      \begin{eqnarray}
Q(0)=Q_{0},~Q'(0)=0,
 \label{eq:oddre}
 \end{eqnarray}
 where       $2<N<4 $ is 
   spatial 
  dimension,   $Q_{0}$ is real,
       \begin{eqnarray}
  Q(\infty) =0,
 \label{eq:oddddddddre}
 \end{eqnarray}
 and energy is zero, i. e.
        \begin{eqnarray}
 H(Q)=\int_{0}^{\infty}\eta^{N-1}\left ( |Q'(\eta)|^{2}-{\frac {1}{2}}  |Q(\eta)|^{4}   
       \right )d\eta  =0.
 \label{eq:odre52c}
 \end{eqnarray} 
 Here, $Q_{0}$ and $a>0$ are constants. Zakharov~\cite{z} and 
 Hastings and McLeod~\cite{has} point out that, because of symmetry and the fact that 
  (\ref{eq:odreaa}) is invariant under a rotation $Q \Rightarrow
 Qe^{i\theta},$  it is justified to assume that     $Q_{0}$ is real and positive.

\medskip \noindent {\bf   Previous Results And Predictions.} In 1988 
  LeMesurier et  al~\cite{L1,L} investigated behavior of solutions when $N=3.$ 
  They analyzed  the     rate at which    
  $Q(\zeta) \to 0$ as $\zeta \to \infty$  for solutions  
  satisfying  (\ref{eq:odreaa})-(\ref{eq:oddre})-(\ref{eq:oddddddddre})-(\ref{eq:odre52c}).
Linearizing (\ref{eq:odreaa}) around the constant solution   $Q=0$ gives 
  \begin{eqnarray}
Q''+{\frac {2}{\zeta}}Q'-Q+ia(Q +\zeta Q')=0,~~0<\zeta<\infty.
 \label{eq:odreaartr1}
 \end{eqnarray}
 They show  that (\ref{eq:odreaartr1}) has independent solutions 
 $Q_{1}$ and $Q_{2}$ which satisfy 
      \begin{eqnarray}
Q_{1} \sim \zeta^{-1- {\frac {i}{a}} }~~{\rm and}~~
Q_{2}\sim \zeta^{ -2+ {\frac {i}{a}}      }e^{-ia{\frac {\zeta^{2}}{2}}}
~{\rm as}~\zeta \to \infty.
 \label{eq:odre52}
 \end{eqnarray}
   LeMesurier et al~\cite{L1,L}  prove that a solution 
    of (\ref{eq:odreaa})-(\ref{eq:oddre}) satisfies   zero-energy condition 
 (\ref{eq:odre52c}) only if, for some $k\ne 0,$  
   \begin{eqnarray}
Q(\zeta) \sim k Q_{1}(\zeta)
~~{\rm as}~~\zeta \to \infty.
\label{eq:oddrde114}
 \end{eqnarray}
 Their numerical   investigation  led to 
  
  \smallskip \noindent  
  {\bf Prediction (I)~\cite{L1,L, sul}} When $N=3$ a wide range of initial conditions exist such that 
  the solution of (\ref{eq:odre55})-(\ref{eq:odre55a})
  asymptotoically approaches~(as $t \to T^{-}$) a self-similar solution of (\ref{eq:odreaa})-(\ref{eq:oddre}) 
when $(Q_{0},a) \approx (1.885,.918).$   The
  solution  satisfies $Q(\infty)=0,$ the   zero-energy condition (\ref{eq:odre52c}), 
  and its profile $|Q|$ is monotonically decreasing.

\medskip \noindent In 1990 Wang~\cite{wang} also investigated the behavior of solutions when  $N=3.$ He proved that, for each $a>0 $ and $Q_{0}>0,$
  the  
  solution  of 
(\ref{eq:odreaa})-(\ref{eq:oddre}) satisfies   
 $Q(\infty)=0.$ However, he did not analyze the number of oscillations in the profile,  $|Q|,$   nor did he determine 
 whether any    solution   satisfies the physically important zero-energy condition (\ref{eq:odre52c}).
In 1995 Kopell and Landman~\cite{kop} proved existence of solutions of    
(\ref{eq:odreaa})-(\ref{eq:oddre})-(\ref{eq:oddddddddre})-(\ref{eq:odre52c})
  when $N>2$ is exponentially close to $N=2.$ 
  In their 1999 book, Sulem and Sulem~\cite{sul} gave an extensive summary of  numerical and theoretical results, and   
  described   physical  relevance of solutions of 
   problem (\ref{eq:odreaa})-(\ref{eq:oddre})-(\ref{eq:oddddddddre})-(\ref{eq:odre52c}).
  In  2000 Budd, Chen and Russel~\cite{budd}   
   obtained further    results.  
     Their  numerical study led to

 \smallskip \noindent  {\bf Prediction (II)~\cite{budd}} When     $2<N<4$   there exists a 
    countably  infinite set of multi-bump solutions
     of (\ref{eq:odreaa})-(\ref{eq:oddre})-(\ref{eq:oddddddddre})-(\ref{eq:odre52c}), and  $n_{a},$ the number of  
     oscillations of the profile  $|Q|,$   satisfies $n_{a} \to \infty$ as $ a \to 0^{+}.$ 
   Furthermore, their  numerical experiments   demonstrate  the important role of  multi-bump solutions  
     in the formation of   singularities in  solutions
    of  problem (\ref{eq:odre55})-(\ref{eq:odre55a}).

     \medskip \noindent 
       Budd et  al~\cite{budd}   also 
  derived basic   properties of  solutions of  (\ref{eq:odreaa})-(\ref{eq:oddre}).
  We will make use of three of   their  results. 
   First, they     showed   that, for any $N \in (2,4),$ $Q_{0}>0$ and   $a>0,$  the solution 
   of (\ref{eq:odreaa})-(\ref{eq:oddre}) exists 
    for all $\zeta \ge 0,$ and 
 an $M>0$ exists such that
          \begin{eqnarray}
 |\zeta Q(\zeta)| \le M~~{\rm and}~~|\zeta^{\alpha}Q'(\zeta)| \le M,~~\zeta \ge 0,
 \label{eq:d3e4doddre}
 \end{eqnarray}
 where    $0<\alpha<N-2$ if $2<N<3$, and $\alpha=1$ if $ 3 \le N<4.$ 
 
 \smallskip \noindent Second, they showed that, if $2<N<4, $ then a solution  of (\ref{eq:odreaa})-(\ref{eq:oddre})-(\ref{eq:oddddddddre}) satisfies 
 \begin{eqnarray}
|\zeta Q'+Q|^{2}+{\frac {1}{2}}\zeta^{2}|Q|^{4}-\zeta^{2}|Q|^{2}=|Q(0)|^{2}-\int_{0}^{\zeta}s|Q(s)|^{4}ds,~~\zeta \ge 0.
   \label{m1aerdddq2dd}
  \end{eqnarray} 
 
 \smallskip \noindent Third, when $2<N<4$ they proved that a solution   of (\ref{eq:odreaa})-(\ref{eq:oddre})-(\ref{eq:oddddddddre}) satisfies   
        \begin{eqnarray}
H(Q)=0 \iff {\Bigg |} \zeta Q'+\left ( 1 +{\frac {i}{a}} \right )Q{\Bigg |}  \to 0~~{\rm as}~~\zeta \to \infty.
 \label{eq:d3e4oddre}
 \end{eqnarray}
 Note that the first inequality in (\ref{eq:d3e4doddre})   implies that $Q(\infty)=0$ for all choices of $Q_{0}>0$ and $a>0,$   which is 
 consistent with  the 
 1990 Wang~\cite{wang} result.
 Subsequently, in 2002 and 2003 Rottschafer and Kaper~\cite{rott, rott2} extended the Kopell et al~\cite{kop} and Budd
 et al~\cite{budd}  results by proving
  existence of families of multi-bump solutions 
  of (\ref{eq:odreaa})-(\ref{eq:oddre})-(\ref{eq:oddddddddre})-(\ref{eq:odre52c}) when 
  $N>2$ is algebraically close to $N=2.$

\medskip \noindent {\bf Goals.} 
 In this paper we focus on dimension $N=3$ and 
 analyze   qualitative behavior of
 solutions of   problem (\ref{eq:odreaa})-(\ref{eq:oddre}). Since the 1984-1986 pioneering investigations by Zhakharhov~\cite{z}, McLaughlin et al~\cite{MCL}, 
    LeMesurier et al~\cite{L1, L} and Landman et al~\cite{Land}, two important theoretical     problems  have remained unresolved: 

\smallskip \noindent {\bf Problem I.}  Do  $Q_{0}>0$ and $a>0$  exist  such that the solution of  
 (\ref{eq:odreaa})-(\ref{eq:oddre})-(\ref{eq:oddddddddre}) satisfies zero-energy condition (\ref{eq:odre52c})~?

\smallskip \noindent {\bf Problem II.}  If  $Q_{0}>0$ and $a>0$ exist such that the solution of   
 (\ref{eq:odreaa})-(\ref{eq:oddre})-(\ref{eq:oddddddddre}) satisfies   (\ref{eq:odre52c}),  can we prove 
 the number of oscillations (i.e. bumps)  of   $|Q|$~?

 \medskip \noindent 
 In order to put   previous numerical predictions on a firm theoretical foundation, it is first necessary to determine 
 the maximal  range of  values   $Q_{0}>0$ and $a>0$ such that    zero-energy
 condition (\ref{eq:odre52c}) is satisfied. Thus, our main goal  
 is to resolve {\bf Problem I.} For this we prove   
   \begin{theo} \label{th1aaa}  Let N=3. For each $Q_{0}>0$ and   $a>0$ the solution of 
   (\ref{eq:odreaa})-(\ref{eq:oddre})   satisfies 
\begin{eqnarray}
 Q(\zeta) \ne 0~~\forall \zeta \ge 0~~{\rm and}~~H(Q)=0.
  \label{m1aderdddq2}
  \end{eqnarray}  
  \end{theo}

\bigskip   \noindent {\bf Discussion}

\smallskip \noindent {\bf (1)} 
 Previous numerical computations~\cite{budd, Land, L1, L, MCL,    z}
    were done on
 finite intervals, hence it is not clear that 
   they   correspond to solutions of  (\ref{eq:odreaa})-(\ref{eq:oddre}) which, for some $Q_{0}>0$ and $a>0,$
   satisfy the physically important    zero-energy condition (\ref{eq:odre52c}) at $\zeta=\infty.$ 
  Theorem~\ref{th1aaa} resolves this important issue since   it is a global result which   guarantees   
 that  zero-energy condition (\ref{eq:odre52c}) holds  
 for all choices $Q_{0}>0$ and $a>0.$ The next step in putting 
    numerical predictions 
 on a firm theoretical footing is to investigate {\bf Problem II} and precisely determine 
    shapes of   profiles of solutions of 
   (\ref{eq:odreaa})-(\ref{eq:oddre})-(\ref{eq:oddddddddre})-(\ref{eq:odre52c}) as $Q_{0}>0$ and $a>0$ vary.
   This will be the object of  future studies.

\smallskip \noindent {\bf (2)} 
The proof of Theorem~\ref{th1aaa} is given in Section~~\ref{sec:q2}.


  \section{Proof of Theorem~\ref{th1aaa}} \label{sec:q2}
 \setcounter{equation}{0} 
 The first step of our   proof of Theorem~\ref{th1aaa} 
  is to  put problem (\ref{eq:odreaa})-(\ref{eq:oddre}) into polar form.
 For this  we    
   substitute  $Q=\rho e^{i\theta}$  into (\ref{eq:odreaa})-(\ref{eq:oddre}), 
  and obtain 
  \begin{eqnarray}
  \rho''+{\frac {2}{\zeta}}\rho'=\rho \left ( (\theta')^{2}+a\zeta\theta'+1 -\rho^{2}\right ), 
     \label{m1aerddd1} \vspace{0.1in} \\
  \theta''+{\frac {2}{\zeta}}  \theta'+2 {\frac {\rho'}{\rho}}\theta' =-a {\frac {(\zeta \rho)'}{\rho}}, 
 \label{m1aerddd2} \vspace{0.1in} \\
  \rho(0)=\rho_{0}>0, ~\rho'(0)=0, ~\theta (0)=0, ~\theta'(0)=0,
  \label{m1aerddd3} 
  \end{eqnarray}
   where  $\rho_{0}>0$ and $a>0.$  
  Combining the zero energy criterion (\ref{eq:d3e4oddre}) with   (\ref{m1aderdddq2}) and the fact that $Q=\rho e^{i\theta},$ 
   we conclude that the proof of Theorem~\ref{th1aaa}  is complete if we show that
   \begin{eqnarray} 
  \rho(\zeta) > 0~~\forall \zeta \ge 0,~~\lim_{\zeta \to \infty}\rho (\zeta)=0,
  \label{m1ade2a}
  \end{eqnarray} 
  and
            \begin{eqnarray}
\lim_{\zeta \to 0} \left ( ( \left ( \zeta \rho \right )')^{2}   +     \left (\zeta \rho \theta' \right )^{2}+
 {\frac {2\zeta \rho^{2} \theta'}{a}}  
+{\frac {\rho^{2}}{a^{2}}}  
  \right )=0.
 \label{eq:d3e356}
 \end{eqnarray}

 \medskip \noindent Our goal in the remainder of this section is to prove that  properties (\ref{m1ade2a})-(\ref{eq:d3e356}) hold.
  For this  we develop seven auxiliary Lemmas in which we prove key qualitative properties of solutions which  
   will allow us to 
    prove    (\ref{m1ade2a})-(\ref{eq:d3e356}). In particular, these technical results  show that, for each $\rho_{0}>0$ and $a>0,$ 
  there exists
  $M>0$ such that the    solution  of 
  initial value probem (\ref{m1aerddd1} )-(\ref{m1aerddd2})-(\ref{m1aerddd3}) satisfies
           \begin{eqnarray}
 0<\zeta \rho  \le M~~{\rm and }~~ 0 \le |\zeta \rho \theta'| \le M~~\forall \zeta \ge 0,~~
  \lim_{\zeta\to \infty} (\zeta\rho)'=0~~{\rm and}~~\lim_{\zeta\to \infty} \theta'=0.
 \label{eq:d22}
 \end{eqnarray} 
  It is easily verified that properties (\ref{m1ade2a})-(\ref{eq:d3e356}) follow   from (\ref{eq:d22}).

\smallskip \noindent   
  
  \begin{lem}  \label{lem2}   Let $\rho_{0}>0$ and $a>0.$  There is an $M>0$ such that the  
     solution of (\ref{m1aerddd1})-(\ref{m1aerddd2})-(\ref{m1aerddd3}) satisfies 
     \begin{eqnarray}
  \label{m123jww1}
  0 \le \zeta \rho \le M,~0 \le    |\zeta \rho'| \le M ~{\rm and}~~0 \le   |\zeta \rho \theta'| \le  M~~\forall \zeta \ge 0.
    \end{eqnarray} \end{lem}

   \noindent  {\bf Proof.}  It follows from  (\ref{eq:d3e4doddre}), combined with the fact that $Q=\rho e^{i\theta},$  that
      \begin{eqnarray}
  \label{m123jww2}
  0 \le \zeta\rho \le M~~{\rm and}~~0 \le \left (\zeta \rho' \right )^{2}+\left (\zeta  \rho \theta' \right )^{2}\le M^{2}
  ~~\forall \zeta \ge 0.
    \end{eqnarray}
 Property (\ref{m123jww1}) follows immediately  from (\ref{m123jww2}). This completes the proof.
 
    \smallskip \noindent {\bf Remark.} The second property in (\ref{m1ade2a})  follows from first property in (\ref{m123jww2}).

  \begin{lem} \label{lem1} 
  Let $\rho_{0}>0$ and $a>0.$ 
  Then the   solution of (\ref{m1aerddd1})-(\ref{m1aerddd2})-(\ref{m1aerddd3}) satisfies 
     \begin{eqnarray}
  \rho(\zeta)>0~~\forall \zeta \ge 0.
  \label{m1a33ee}
  \end{eqnarray}
\end{lem}
   \noindent {\bf Remark.} It follows from (\ref{m1a33ee}) that $|Q(\zeta)|=\rho>0~~\forall \zeta \ge 0,$ hence 
     $Q(\zeta) \ne 0~~\forall \zeta \ge 0.$
This  proves the first property in  (\ref{m1ade2a}), and also the first property in  (\ref{m1aderdddq2}).

\smallskip \noindent 
 {\sl Proof.} Suppose, for contradiction,  that  $\bar \zeta>0$
 exists such that
   \begin{eqnarray}
   \rho(\zeta)>0~~\forall \zeta \in [0,\bar \zeta )~~{\rm and }~~ \rho ( \bar \zeta)=0.
  \label{m144dddee}
   \end{eqnarray}
The first step in obtaining a contradiction   to (\ref{m144dddee}) is to recall  from (\ref{m123jww1}) that 
  \begin{eqnarray}
  \label{m123jdww1}
  |\zeta \rho \theta'| \le  M~~\forall \zeta \in [0,\bar \zeta].
    \end{eqnarray}
  Next, we  write (\ref{m1aerddd2}) as
   \begin{eqnarray}
     \left ( \zeta^{2}\rho^{2}\theta' \right )'= -{\frac {a\zeta}{2}} \left ( \left ( \zeta \rho\right )^{2}\right )'.
  \label{m1ae345}
  \end{eqnarray}
  Integrating (\ref{m1ae345}) by parts gives 
        \begin{eqnarray}
      \zeta^{2}\rho^{2} \left ( \theta'  +{\frac {a\zeta}{2}} \right )= {\frac {a}{2}} \int_{0}^{\zeta}\left
       (t \rho(t)\right )^{2}dt. 
  \label{m1aer346}
  \end{eqnarray}
  Now define the   finite, positive  value
           \begin{eqnarray}
       C= {\frac {a}{2  {\bar \zeta}^{2}}} \int_{0}^{\bar \zeta}\left
       (t \rho(t)\right )^{2}dt>0.
   \label{m133r346}
   \end{eqnarray}
  We conclude from (\ref{m1aer346}), (\ref{m133r346}) and the fact that $\rho(\zeta) \to 0^{+}$ as $\zeta \to {\bar \zeta}^{-},$ that
           \begin{eqnarray}
      \theta'\sim {\frac {C}{\rho^{2}}} ~~{\rm as}~~\zeta \to {\bar \zeta}^{-}.
   \label{m13d357}
  \end{eqnarray}
  It follows from (\ref{m13d357}) that $|\zeta \rho \theta'| \to \infty$ as $\zeta \to {\bar \zeta}^{-},$ which 
  contradicts (\ref{m123jdww1}). This completes the proof of Lemma~\ref{lem1}.

       \begin{lem}  \label{lem4}   
      Let $\rho_{0}>0$ and $a>0.$    Then the     solution of (\ref{m1aerddd1})-(\ref{m1aerddd2})-(\ref{m1aerddd3}) satisfies 
   \begin{eqnarray}
 \int_{0}^{\zeta} \left ( t\rho(t) \right )^{2}dt \to \infty~~{\rm as}~~ \zeta \to \infty.
 \label{m1ad12aa1}
  \end{eqnarray}
 \end{lem}

 \bigskip \noindent {\bf Proof.}  We assume, for contradiction, that   $\bar \rho_{0}>0$ and $\bar a>0$ exists such that the 
 solution of (\ref{m1aerddd1})-(\ref{m1aerddd2})-(\ref{m1aerddd3}) corresponding to $(\rho_{0},a)=(\bar \rho_{0},\bar a)$ satisfies
    \begin{eqnarray}
0<D= \int_{0}^{\infty} \left ( t\rho(t) \right )^{2}dt <\infty.
 \label{m1ad12ffd}
  \end{eqnarray}
 The first step in obtaining a contradiction to   (\ref{m1ad12ffd}) is to write   equation (\ref{m1aer346}) as
  \begin{eqnarray}
      \zeta \rho   \left  (  \zeta \rho \theta'\right )  +{\frac {a}{2}}\zeta \left (\zeta \rho \right )^{2}
        = {\frac {a}{2}} \int_{0}^{\zeta}\left
       (t \rho(t) \right )^{2}dt. 
  \label{m1aer346yy}
  \end{eqnarray}
 It follows from (\ref{m123jww1}) that    $M>0$ exists such that the first term in (\ref{m1aer346yy})     satisfies
  \begin{eqnarray}
  \label{m123jwdw}
  0\le  |\zeta \rho   \left  (  \zeta \rho \theta'\right )| \le M^{2}~~\forall \zeta \ge0.
  \end{eqnarray}
 Next, we claim that  
   \begin{eqnarray}
  \label{m123jwdddw}
  \zeta \rho (\zeta) \to 0~~{\rm as}~~\zeta \to \infty.  
  \end{eqnarray}
  If (\ref{m123jwdddw}) is false, there exist  $\delta \in (0,M)$ and a positive, increasing, unbounded sequence 
  $\left ( \zeta_{N} \right )$ exists such that
     \begin{eqnarray}
  \label{m123jwdd}
  0<\delta < \zeta_{N} \rho (\zeta_{N}) \le M~~\forall N \ge 1.
  \end{eqnarray}
  Along the sequence  $\left ( \zeta_{N} \right )$ equation (\ref{m1aer346yy}) becomes
   \begin{eqnarray}
      \zeta_{N} \rho (\zeta_{N})  \left  (  \zeta_{N} \rho (\zeta_{N}) \theta'(\zeta_{N})\right )
        +{\frac {a}{2}}\zeta_{N} \left (\zeta_{N} \rho (\zeta_{N}) \right )^{2}
        = {\frac {a}{2}} \int_{0}^{\zeta_{N}}\left
       (t \rho(t) \right )^{2}dt. 
  \label{m1a999}
  \end{eqnarray}
  Combining    the   bounds    in (\ref{m123jww1}) and 
   (\ref{m123jwdd}) with the left side of (\ref{m1a999}) gives
     \begin{eqnarray}
      \zeta_{N} \rho (\zeta_{N})  \left  (  \zeta_{N} \rho (\zeta_{N}) \theta'(\zeta_{N})\right )
        +{\frac {a}{2}}\zeta_{N} \left (\zeta_{N} \rho (\zeta_{N}) \right )^{2}
       \ge  {\frac {a}{2}}\zeta_{N} \delta^{2}-M^{2}~~\forall N \ge 1.
  \label{m1a9d99dd}
  \end{eqnarray}
  It follows from (\ref{m1a9d99dd}) and the fact that $\zeta_{N} \to \infty$ as $N \to \infty$ that
      \begin{eqnarray}
     \lim_{N \to \infty} \left ( \zeta_{N} \rho (\zeta_{N})  \left  (  \zeta_{N} \rho (\zeta_{N}) \theta'(\zeta_{N})\right )
        +{\frac {a}{2}}\zeta_{N} \left (\zeta_{N} \rho (\zeta_{N}) \right )^{2} \right )=\infty.
   \label{m1a9432}
  \end{eqnarray}
  However, the right side of (\ref{m1a9d99dd}) remains bounded as $N \to \infty,$ since, by (\ref{m1ad12ffd}),  
   \begin{eqnarray}
     0< \lim_{N \to \infty} \left ( {\frac {a}{2}} \int_{0}^{\zeta_{N}}\left
       (t \rho(t) \right )^{2}dt \right )={\frac {aD}{2}} <\infty. 
  \label{m1a999a}
  \end{eqnarray}
 We conclude from (\ref{m1a999}), (\ref{m1a9432}) and (\ref{m1a999a})   that the left side of 
  (\ref{m1a999}) becomes unbounded as $N \to \infty,$ whereas the right side remains bounded as $N \to \infty,$ a contradiction. 
  Thus,  property    (\ref{m123jwdddw}) holds, as claimed.
  It now follows from (\ref{m123jwdddw}), and the fact that $0 \le |\zeta \rho \theta'| \le M~~\forall \zeta \ge 0$~(i.e. see  
  (\ref{m123jww1})),  that the first term in 
       (\ref{m1aer346yy}) satisfies 
        \begin{eqnarray}
     \zeta \rho   \left  (  \zeta \rho \theta'\right )  \to 0~~{\rm as}~~ \zeta \to \infty.
  \label{m1a99d9aggg}
  \end{eqnarray}
  We conclude from (\ref{m1ad12ffd}), (\ref{m1aer346yy}) and (\ref{m1a99d9aggg}) that
         \begin{eqnarray}
  \zeta \left (\zeta \rho \right )^{2}  \to D~~{\rm as}~~ \zeta \to \infty.
  \label{m1a99d9afff}
  \end{eqnarray}
 It follows from (\ref{m1a99d9afff}) and the   assumption   $D>0,$ that $\zeta_{1}>0,$  $k_{1}>0$ exist such that
           \begin{eqnarray}
   \left (\zeta \rho \right )^{2}  \ge {\frac {k_{1}}{  \zeta}}~~\forall \zeta \ge \zeta_{1}.
  \label{m1arrr}
  \end{eqnarray}
 We conclude from (\ref{m1arrr}) that
            \begin{eqnarray}
	    \int_{0}^{\zeta} \left ( t\rho(t) \right )^{2}dt \ge  \int_{\zeta_{1}}^{\zeta} {\frac {k_{1}}{t}}dt = k_{1} 
	  \left (   \ln(\zeta)-\ln(\zeta_{2}) \right )~~\forall \zeta \ge \zeta_{1}. 
   \label{m1arrra}
  \end{eqnarray}
 It follows from (\ref{m1arrra}) that $\int_{0}^{\zeta} \left ( t\rho(t) \right )^{2}dt \to \infty$ as $\zeta \to \infty,$   contradicting
   (\ref{m1ad12ffd}). Thus, we conclude that   (\ref{m1ad12aa1}) 
    holds as claimed. This completes the proof of Lemma~\ref{lem4}.

       \begin{lem}  \label{lem6}   
      Let $\rho_{0}>0$ and $a>0.$  Then the     solution of (\ref{m1aerddd1})-(\ref{m1aerddd2})-(\ref{m1aerddd3}) satisfies 
    \begin{eqnarray}
 \lim_{\zeta \to \infty} \zeta \left (\zeta \rho \right )^{2} =\infty,~~\lim_{\zeta \to \infty}
   {\frac { \left ( \zeta \rho \right )' }{\rho}} = 0~~{\rm and}~~\lim_{\zeta \to \infty} \left ( \zeta \rho \right )' =0.
 \label{m1ad1tgggg}
  \end{eqnarray}  
     \end{lem}
    {\bf Remark.} The first two limits   in (\ref{m1ad1tgggg}) will be used to prove the  third limit, 
    and  third limit proves the fourth property in (\ref{eq:d22}).
     Properties (\ref{m1ad1tgggg}) will   play an essential 
     role in completing the proof that  $\lim_{\zeta \to \infty}\theta'=0,$ the crucial fifth 
    property in (\ref{eq:d22}).

 \medskip \noindent {\bf Proof.} First, it follows from 
 (\ref{m1ad12aa1}) and (\ref{m1aer346yy}) that 
  \begin{eqnarray}
  \lim_{\zeta \to \infty} \left (    \zeta \rho   \left  (  \zeta \rho \theta'\right )  +{\frac {a}{2}}\zeta \left (\zeta \rho \right )^{2} \right )
        =  \lim_{\zeta \to \infty} \left ( {\frac {a}{2}} \int_{0}^{\zeta}\left
       (t \rho(t) \right )^{2}dt \right ) = \infty. 
  \label{m1ae76}
  \end{eqnarray}
Recall from (\ref{m123jwdw}) that $ 0\le  |\zeta \rho   \left  (  \zeta \rho \theta'\right )| \le M^{2}~~\forall \zeta \ge 0.$ 
This     and (\ref{m1ae76}) imply   that 
       \begin{eqnarray}
 \lim_{\zeta \to \infty} \zeta \left (\zeta \rho \right )^{2} =\infty,
  \label{m1ae76a}
  \end{eqnarray}
which proves the first property in (\ref{m1ad1tgggg}).
 Next, divide   (\ref{m1aer346yy}) by ${\frac {a}{2}}\zeta\left (\zeta \rho \right )^{2}$ and get
  \begin{eqnarray}
     {\frac {2}{a}} {\frac {\zeta \rho   \left  (  \zeta \rho \theta'\right )}{  \left (\zeta \left (\zeta \rho \right )^{2} \right )}}  +
     1   =  {\frac { \int_{0}^{\zeta}\left
       (t \rho(t) \right )^{2}dt}{\left (\zeta \left (\zeta \rho \right )^{2} \right )}},~~\zeta >0. 
  \label{m1a555}
  \end{eqnarray}
It follows from (\ref{m1ae76a}) and the bound $ 0\le  |\zeta \rho   \left  (  \zeta \rho \theta'\right )| \le M^{2}~~\forall \zeta \ge 0$ that
  \begin{eqnarray}
    \lim_{\zeta \to \infty} {\frac {2}{a}} {\frac {\zeta \rho  
     \left  (  \zeta \rho \theta'\right )}{  \left (\zeta \left (\zeta \rho \right )^{2} \right )}}   
      =0.
  \label{m1a5rre}
  \end{eqnarray}
Taking the limit as $\zeta \to \infty$ to both sides of (\ref{m1a555}), and   using (\ref{m1a5rre}), gives
  \begin{eqnarray}
     1   = \lim_{\zeta \to \infty} {\frac { \int_{0}^{\zeta}\left
       (t \rho(t) \right )^{2}dt}{\left (\zeta \left (\zeta \rho \right )^{2} \right )}}. 
  \label{m1da555}
  \end{eqnarray}
It follows from (\ref{m1ad12aa1}),   (\ref{m1ae76a}), (\ref{m1da555})   
  and L'Hopital's rule  that
      \begin{eqnarray}
     1 = 
     \lim_{\zeta \to \infty}  {\frac {\left (\zeta\rho(\zeta) \right )^{2}}{ \left (\zeta\rho(\zeta) \right )^{2}+
     2\zeta \left (\zeta\rho(\zeta) \right )\left (\zeta\rho(\zeta) \right )'}} =
       \lim_{\zeta \to \infty}  {\frac {1}{ 1+
     2 \left (\zeta\rho(\zeta) \right )'/\rho}},
   \label{m1ade}
  \end{eqnarray} 
hence
        \begin{eqnarray}
  \lim_{\zeta \to \infty}
    {\frac { \left (\zeta\rho(\zeta) \right )'}{\rho}}=0,
   \label{m1add54}
  \end{eqnarray}
which proves the second property in (\ref{m1ad1tgggg}). Finally, we conclude from (\ref{m1add54})  and  the bound $0 \le \zeta \rho \le M~~\forall \zeta \ge 0$
   that
 $\zeta^{*}>0$ exists such that
    \begin{eqnarray}
   |\left (\zeta\rho(\zeta) \right )'| \le \rho \le {\frac {M}{\zeta}}~~\forall \zeta \ge \zeta^{*}. 
  \label{m1dd54}
  \end{eqnarray}
It follows from (\ref{m1dd54}) that $\left (\zeta\rho(\zeta) \right )' \to 0$ as $\zeta \to \infty,$ which proves the  
third property in (\ref{m1ad1tgggg}). This completes the proof of   Lemma~\ref{lem6}.

 
 \noindent The remainder of this section is devoted to proving the fifth  property in 
  (\ref{eq:d22}), i.e. 
  for every $\rho_{0}>0$ and $a>0,$ the   
  solution of  (\ref{m1aerddd1})-(\ref{m1aerddd2})-(\ref{m1aerddd3}) satisfies 
  $\lim_{\zeta \to \infty} \theta' =0.$  In the next  three Lemmas we    prove this property   
    by using 
  the results proved above, and  we   also make    extensive use of the functional 
     \begin{eqnarray}
 E= \left ( \left ( \zeta \rho  \right )' \right )^{2}+ 
  \left ( \zeta \rho  \right )^{2} \left ( (\theta')^{2}-1 \right ) +{\frac {1}{2}}\zeta^{2}\rho^{4},
  \label{m1aqq}
  \end{eqnarray}
which satisfies, because  $\rho>0~~\forall \zeta>0, $
\begin{eqnarray}
 E(0)=\rho_{0}^{2}~~{\rm and}~~ E'= -\zeta\rho^{4}<0~~\forall \zeta > 0.
 \label{m1adedtryaac}
  \end{eqnarray}
 An integration gives 
       \begin{eqnarray}
  \left ( \left ( \zeta \rho  \right )' \right )^{2}+ 
  \left ( \zeta \rho  \right )^{2} \left ( (\theta')^{2}-1 \right ) +{\frac
  {1}{2}}\zeta^{2}\rho^{4}=\rho_{0}^{2}-\int_{0}^{\zeta}t\rho(t)^{4}dt~~\forall \zeta \ge 0.
  \label{m1adedddtryaaj}
  \end{eqnarray}
 {\bf Remark.} Equation~(\ref{m1adedddtryaaj}) is the same
  as  equation~(\ref{m1aerdddq2dd}) derived    by Budd et
 al~\cite{budd}.

  \medskip \noindent 
  The  key properties of the functional $E$  
  are proved in the next
  two technical Lemmas.

    \begin{lem}  \label{lem3}   
    Let $\rho_{0}>0$ and $a>0.$   Then the     solution of (\ref{m1aerddd1})-(\ref{m1aerddd2})-(\ref{m1aerddd3}) satisfies
   \begin{eqnarray}
 -\infty < E(\infty) \le 0.
 \label{m1ad12}
  \end{eqnarray}
     \end{lem}
 {\bf Proof.}
 We  conclude from  (\ref{m1adedtryaac}), Lemma~\ref{lem1} 
 and the bound  $0<\zeta\rho \le M~\forall \zeta \ge 0$     that
                 \begin{eqnarray}
    E(\infty) = \rho_{0}^{2}-\int_{0}^{\infty}t\rho^{4}(t)dt>-\infty.
  \label{m55dde}
  \end{eqnarray}
It remains to prove  that   
  \begin{eqnarray}
     E(\infty) \le 0.
 \label{m1ad1d2ww456}
  \end{eqnarray}
  Suppose, for contradiction, 
  that $\bar \rho_{0}>0$ and $\bar a>0$ exist such that, when $ (\rho_{0},a)=(\bar \rho_{0},\bar a),$   the 
 solution of   (\ref{m1aerddd1})-(\ref{m1aerddd2})-(\ref{m1aerddd3}) satisfies 
    \begin{eqnarray}
  E(\infty) =\bar \lambda > 0.
 \label{m1ad1d2wwwhj}
  \end{eqnarray}
  The first step in obtaining a contradiction to  (\ref{m1ad1d2wwwhj}) is 
 to write   (\ref{m1aqq}) as
  \begin{eqnarray}
  \left ( \zeta \rho  \right )^{2} \left ( (\theta')^{2}-1 \right )  
  = E(\zeta)-\left ( \left ( \zeta \rho  \right )' \right )^{2} -{\frac {1}{2}}\zeta^{2}\rho^{4}, ~~\zeta \ge 0.
  \label{m1adedtryaajwh}
  \end{eqnarray}
 From the third property in (\ref{m1ad1tgggg})   and the fact that    $0\le \zeta\rho\le M $ 
 we conclude that 
  \begin{eqnarray}
 \lim_{\zeta \to \infty} \left (\left ( \zeta \rho \right )'\right )^{2} =0~~{\rm and}~~\lim_{\zeta \to \infty} {\frac {1}{2}}\zeta^{2}\rho^{4}  =0.
 \label{m1adret}
  \end{eqnarray} 
 It follows from   (\ref{m1ad1d2wwwhj}),   (\ref{m1adedtryaajwh}) and (\ref{m1adret})
that   
  $\lim_{\zeta \to \infty}\left (\left ( \zeta \rho  \right )^{2} \left ( (\theta')^{2}-1 \right )  \right ) = {\bar \lambda}>0.$
 Thus, there exists $\bar \zeta>0$ such that 
    \begin{eqnarray}
   \left ( \zeta \rho  \right )^{2} \left ( (\theta')^{2}-1 \right )  
  \ge {\frac {{\bar \lambda}}{2}}~~\forall \zeta \ge  {\bar \zeta}. 
  \label{m1aded983a}
  \end{eqnarray}
 We conclude from (\ref{m1aded983a}) that  
 $   \left ( \zeta \rho  \right )^{2} (\theta')^{2}  \ge {\frac {{\bar \lambda}}{2  }} ~~\forall \zeta \ge  {\bar \zeta}.$
Since $\zeta\rho \theta'$ is continuous,  either
        \begin{eqnarray}
 \zeta\rho \theta'\ge \left ({\frac {{\bar \lambda}}{2  }}\right )^{1/2}~~\forall \zeta \ge  {\bar \zeta}~~{\rm or}~~\zeta \rho \theta'\le-
\left ( {\frac {{\bar \lambda}}{2  }}\right )^{1/2}~~\forall \zeta \ge  {\bar \zeta}.
  \label{m1art5rt}
  \end{eqnarray}
  Also, since $0 \le \zeta \rho \le M~~\forall \zeta \ge 0,$ then    $\bar \zeta$ can be chosen large enough  so that 
       \begin{eqnarray}
 1-\rho^{2}>0~~\forall \zeta \ge {\bar \zeta}.
  \label{m1artd}
  \end{eqnarray}
 Now  write (\ref{m1aerddd2}) as 
  \begin{eqnarray}
  \left (\zeta \rho\right )'' =\zeta \rho \left ( (\theta')^{2}+a\zeta\theta'+1 -\rho^{2}\right ).
     \label{m1a57}  
  \end{eqnarray}
 Suppose   that $   \zeta\rho \theta'\ge \left ({\frac {{\bar \lambda}}{2  }}\right )^{1/2}~~\forall \zeta \ge  {\bar \zeta}.$
Then it  follows from   (\ref{m1artd}) and (\ref{m1a57})  that
  \begin{eqnarray}
  \left (\zeta \rho\right )'' \ge     a \zeta  \left ({\frac {{\bar \lambda}}{2  }}\right )^{1/2}  ~~\forall \zeta \ge  {\bar \zeta}. 
     \label{m1a5d7}  
  \end{eqnarray}
An integration from $\bar \zeta$ to $\zeta$ gives
  \begin{eqnarray}
   \left (\zeta \rho\right )' \ge  \left (\zeta \rho\right )'|_{\zeta=\bar \zeta}+
   a   \left ({\frac {{\bar \lambda}}{2  }}\right )^{1/2} 
    \left (  {\frac { \zeta^{2}}{2}}- {\frac {  {\bar \zeta}^{2}}{2}} \right )~~\forall \zeta \ge  {\bar \zeta}. 
     \label{m1a5d7dd}  
  \end{eqnarray}
It follows from (\ref{m1a5d7dd}) that $\left (\zeta \rho\right )' \to \infty$ as $\zeta \to \infty,$ contradicting the third property in 
(\ref{m1ad1tgggg}). Thus, the first possibility in (\ref{m1art5rt}) cannot occur.
 It remains to assume, for contradiction,  that  the second possibility in (\ref{m1art5rt}) occurs, i.e. 
        \begin{eqnarray}
\zeta \rho \theta'\le-
\left ( {\frac {{\bar \lambda}}{2  }}\right )^{1/2}~~\forall \zeta \ge  {\bar \zeta}.
  \label{m143t}
  \end{eqnarray}
 First, we conclude   from (\ref{m1a57}) and the bound  $0\le \zeta \rho \le M~~\forall \zeta \ge 0,$ that 
  \begin{eqnarray}
  \left (\zeta \rho\right )'' \le \zeta \rho \theta' \left (  \theta'  +a\zeta  \right )+\zeta\rho \le 
  \zeta \rho \theta' \left (  \theta' +a\zeta  \right )+M~~\forall \zeta \ge 0.
     \label{m1a5d7b}  
  \end{eqnarray}
  Also,    from (\ref{m1aer346}) it follows that 
           \begin{eqnarray}
   \theta' + {\frac {a\zeta}{2}} \ge 0~~\forall \zeta \ge 0.
    \label{m1aeez}
   \end{eqnarray}
 Combining this property with (\ref{m143t}) and   (\ref{m1a5d7b}), we obtain 
   \begin{eqnarray}
  \left (\zeta \rho\right )''   \le 
  \zeta \rho \theta' \left (  \theta'  +a\zeta  \right )+
  M\le-\left ( {\frac {{\bar \lambda}}{2  }}\right )^{1/2}{\frac {a\zeta}{2}}+M ~~\forall \zeta \ge {\bar \zeta}.
     \label{m1a5dd7}  
  \end{eqnarray}
Integrating from $\bar \zeta$ to $\zeta,$ we conclude that 
$ \left (\zeta \rho\right )' \to -\infty$ as $\zeta \to \infty,$   contradicting the   third property in  (\ref{m1ad1tgggg}).
Since     suppostion  (\ref{m1ad1d2wwwhj}) has led  to a contradiction, we conclude that 
(\ref{m1ad1d2wwwhj}) cannot occur, hence   $E(\infty) \le 0.$  This completes the proof of 
Lemma~\ref{lem3}.

 \medskip \noindent In order to complete the proof of Theorem~\ref{th1aaa}   we first need to eliminate the possibility that  
   a solution of (\ref{m1aerddd1})-(\ref{m1aerddd2})-(\ref{m1aerddd3}) satisfies 
  $E(\infty)=0.$ We do this in

    \begin{lem}  \label{lem3a}   
    Let $\rho_{0}>0$ and $a>0.$   Then the     solution of (\ref{m1aerddd1})-(\ref{m1aerddd2})-(\ref{m1aerddd3}) satisfies 
   \begin{eqnarray}
  E(\infty)<0.
 \label{m1ad12aab}
  \end{eqnarray}
     \end{lem}
 {\bf Proof.}
   Suppose, for contradiction,  that $\bar \rho_{0}>0,$   $\bar a>0$ exist such that     
  property (\ref{m1ad12aab})
  does not hold when $ (\rho_{0},a)=(\bar \rho_{0},\bar a).$  This supposition and  Lemma~\ref{lem3} imply that
     \begin{eqnarray}
  E(\infty)=0.
 \label{m1ad12aa}
  \end{eqnarray}
  The first step in obtaining a contradiction to (\ref{m1ad12aa}) is to show that
       \begin{eqnarray}
  \lim_{\zeta \to \infty} \left (\zeta\rho\theta'\right )^{2}=0.
   \label{m1rrfd2aa}
  \end{eqnarray}
 Before proving property (\ref{m1rrfd2aa}) we show how it leads to a contradiction of (\ref{m1ad12aa}). First, recall that the   
   third property in  (\ref{m1ad1tgggg}) is 
      \begin{eqnarray}
     \lim_{\zeta \to \infty} \left ( \left  ( \zeta \rho \right )' \right )^{2} =0.
   \label{m1rrdd}
  \end{eqnarray} 
  Substituting (\ref{m1rrfd2aa}) and (\ref{m1rrdd}) into the left side of   equation  (\ref{eq:d3e356}), and using the bound $0 \le
  \zeta \rho\le M~~\forall \zeta \ge 0,$ we conclude that    zero-energy conditon (\ref{eq:d3e356}) is satisfied, i.e.  
     \begin{eqnarray}  
\lim_{\zeta \to 0} \left (  \left ( \zeta \rho \right )')^{2}   + \left ( \rho^{2} \left (\zeta^{2}(\theta')^{2}+{\frac {2\zeta\theta'}{a}}
+{\frac {1}{a^{2}}} \right ) \right
)^{2} \right )^{\frac {1}{2}} =0.
 \label{eq:d3e356aab}
 \end{eqnarray} 
 It follows from (\ref{eq:d3e4oddre}) and (\ref{eq:d3e356aab}) that
      $H(Q)=0$ where energy 
   $H(Q)$ is defined  in (\ref{eq:odre52c}).
     In turn, as we pointed out  in Section~\ref{sec:intro} 
     LeMesurier et al~\cite{L1,L} proved that    when  
      $H(Q)=0$
   there exists  $k\ne 0$    such that   
      \begin{eqnarray}
  Q(\zeta) \sim k \zeta^{-1} e^{-i\ln(\zeta)/a}~~{\rm as }~~\zeta \to \infty.
  \label{eq:oddddrde114a}
  \end{eqnarray}
 From (\ref{eq:oddddrde114a}),  and the fact that $Q=\rho e^{i\theta},$  it follows that 
     \begin{eqnarray}
 |\zeta Q|=\zeta\rho \to |k|~~{\rm as}~~\zeta \to \infty.
  \label{eq:oddddda}
  \end{eqnarray}
  Substituting    (\ref{m1rrfd2aa}), (\ref{m1rrdd})  and (\ref{eq:oddddda})
   into the left side of (\ref{m1adedtryaajwh}), we     conclude that
    \begin{eqnarray}
   \lim_{\zeta \to \infty} \left ( \left ( \left ( \zeta \rho  \right )' \right )^{2} +
 \left ( \left ( \zeta \rho  \right )^{2} \left ( (\theta')^{2}-1 \right ) \right )\right ) =-|k|^{2}.
 \label{m1adedw}
  \end{eqnarray}
 However, because of (\ref{m1ad12aa}) and properties (\ref{m1adret}) and  (\ref{m1ad12aa}),   the right side of  (\ref{m1adedtryaajwh}) 
  tends to zero as $\zeta \to \infty,$ a contradiction.
   Thus, it  remains to prove  (\ref{m1rrfd2aa}). We assume, for contradiction, that 
property  (\ref{m1rrfd2aa}) does not hold, hence
        \begin{eqnarray}
  \limsup_{\zeta \to \infty} \left (\zeta\rho\theta'\right )^{2}>0.
   \label{m1rrfd2aa12}
  \end{eqnarray}
From  (\ref{m1rrfd2aa12})  
  and the fact that $|\zeta\rho\theta'| \le M~~\forall \zeta \ge 0,$ 
  we conclude that  a value $\delta >0$ exists, and also a positive, increasing, unbounded sequence $ \left (\zeta_{N} \right )$ exists   such that
      \begin{eqnarray}
   4\delta^{2} \le \left ( \zeta_{N}\rho(\zeta_{N})\theta' (\zeta_{N} )\right )^{2}  \le M^{2}~~\forall N \ge 1,
    \label{m1rrfd2}
   \end{eqnarray} 
 and therefore, by considering subsequences if necessary,  either 
       \begin{eqnarray}
  -M \le    \zeta_{N}\rho(\zeta_{N})\theta'(\zeta_{N}) \le -2\delta ~~\forall N \ge 1~~{\rm or}~~2\delta \le 
    \zeta_{N}\rho(\zeta_{N})\theta'(\zeta_{N}) \le M~~\forall N \ge 1
  \label{m1rrfd2ac}
  \end{eqnarray}
  Thus, the proof of Lemma~\ref{lem3a} is complete if we eliminate each of the cases given in (\ref{m1rrfd2ac}). 
  In order to eliminate  these two cases, we first 
 need to derive  upper and lower bounds on   $\zeta_{N}\rho(\zeta_{N})$ when $N\gg 1.$   For this we evaluate 
    (\ref{m1aer346yy}) at $\zeta=\zeta_{N}$ and get 
   \begin{eqnarray}
  \left ( \zeta \rho  \right )^{2} \left ( (\theta')^{2}-1 \right ) { \Big |}_{\zeta=\zeta_{N}} 
  = E(\zeta_{N})-\left ( \left ( \zeta \rho  \right )' \right )^{2}{ \Big |}_{\zeta=\zeta_{N}} -{\frac {1}{2}}\zeta_{N}^{2}\rho(\zeta_{N})^{4}~~\forall N
  \ge 1.
  \label{m1adedtryaajw11}
  \end{eqnarray}
From  (\ref{m1adret}) and  (\ref{m1ad12aa}) it follows that the right side of (\ref{m1adedtryaajw11}) tends to zero as $N \to \infty.$
Thus, the left  side of (\ref{m1adedtryaajw11}) satisfies 
  \begin{eqnarray}
 \lim_{N \to \infty} \left ( \left ( \zeta \rho  \right )^{2} \left ( (\theta')^{2}-1 \right )\right ) {\Big |}_{\zeta=\zeta_{N}} =0.
  \label{m1adedtryaaj11}
   \end{eqnarray}
We conclude from (\ref{m1rrfd2}),   (\ref{m1adedtryaaj11}), and the bound $ 0 \le \zeta\rho\le M~~\forall \zeta \ge 0,$ that
  \begin{eqnarray}
 \delta  \le    \zeta_{N} \rho(\zeta_{N})    \le M ~~{\rm when}~~N \gg 1.
  \label{m1a34}
  \end{eqnarray}
    Also, it follows from the fact that $0 \le \zeta\rho \le M~~\forall \zeta \ge 0$ and  property (\ref{m1ad1tgggg}) in Lemma~\ref{lem6} 
   that ${\hat \zeta}>0$ exists such that
    \begin{eqnarray}
 0<\rho^{2}(\zeta)<1~~\forall \zeta >\hat \zeta,
  \label{m1a34tt} 
  \end{eqnarray}
   \begin{eqnarray}
     -{\frac {\delta}{2}} \le \left ( \zeta \rho \right )'
   \le {\frac {\delta}{2}}~~{\rm and}~~
 -{\frac {\delta^{2}}{aM^{2}}} \le 
   {\frac { \left ( \zeta \rho \right )' }{\rho}} \le {\frac {\delta^{2}}{aM^{2}}}~~\forall \zeta >{\hat \zeta}.
 \label{m1ad1tgggg14}
  \end{eqnarray} 
 We assume, without loss of genrality, that     $\left (\zeta_{N}\right )$ is chosen so that 
 $\zeta_{N} > {\hat \zeta}~~\forall N \ge 1,$
 and also  (\ref{m1a34})  holds  $\forall N \ge 1.$

  \smallskip \noindent We now consider the  two cases in (\ref{m1rrfd2ac}). Suppose,  first of all,  that 
         \begin{eqnarray}
    2\delta \le \zeta_{N}\rho(\zeta_{N})\theta'(\zeta_{N}) \le M~~\forall N \ge 1.
   \label{m1rrfd2b}
  \end{eqnarray} 
  Next, we  derive a lower bound for $\zeta\rho\theta'.$ For this  write (\ref{m1aerddd2}) as
   \begin{eqnarray}
     \left ( \zeta^{2}\rho^{2}\theta' \right )'= -a \left (\zeta \rho \right )^{2}
     {\frac { \left (   \zeta \rho   \right )'}{\rho}}.
  \label{m1art}
  \end{eqnarray} 
  From (\ref{m1art}),   the   bound  ${\frac { \left ( \zeta \rho \right )' }{\rho}} \le {\frac {\delta^{2}}{aM^{2}}}$ 
  in (\ref{m1ad1tgggg14}),   
and the bound  $ 0\le \zeta\rho\le M$ we get
   \begin{eqnarray}
     \left ( \zeta^{2}\rho^{2}\theta' \right )'\ge  - \delta^{2}~~\forall \zeta > \hat \zeta.
  \label{m1artabb}
  \end{eqnarray}
It follows from  (\ref{m1a34}), (\ref{m1rrfd2b})  and an integration of (\ref{m1artabb}) that, for each $N \ge 1,$ 
   \begin{eqnarray}
    \zeta^{2}\rho^{2}\theta'  \ge  \zeta_{N}^{2} \rho^{2} (\zeta_{N}) \theta'(\zeta_{N})   - \delta^{2}(\zeta-\zeta_{N})
    \ge 2\delta^{2}- \delta^{2} = \delta^{2},~~\zeta_{N}\le \zeta \le  \zeta_{N}+1.
  \label{m1arta}
  \end{eqnarray}
Next,  write (\ref{m1aerddd2}) as 
 \begin{eqnarray}
  (\zeta \rho)'' =\zeta \rho \left ( (\theta')^{2}+a\zeta\theta'+1 -\rho^{2}\right ).
     \label{m1a456rr} 
  \end{eqnarray}
 It follows from (\ref{m1a34tt}), (\ref{m1arta}), the fact that $0< \zeta \rho \le M~~\forall \zeta> 0$    and (\ref{m1a456rr})  that 
  solution of (\ref{m1aerddd1})-(\ref{m1aerddd2})-(\ref{m1aerddd3}) satisfies
  \begin{eqnarray}
  (\zeta \rho)'' \ge  a \zeta   {\frac {\left (  \zeta \rho\ \right )^{2} \theta'}{\zeta \rho}} \ge {\frac {a \delta^{2}\zeta}{M}}, ~~
  \zeta_{N} \le \zeta \le \zeta_{N}+1, ~N \ge 1
     \label{m1a456rrb} 
  \end{eqnarray}
An integration of (\ref{m1a456rrb}) from $\zeta_{N}$ to $\zeta$ gives
  \begin{eqnarray}
  (\zeta \rho)'  \ge   (\zeta \rho)'{\Big |}_{\zeta=\zeta_{N}}+  {\frac {a \delta^{2}}{M}} \left ( {\frac {\zeta^{2}}{2}}-
  {\frac {\zeta_{N}^{2}}{2}} \right ),~~ 
  \zeta_{N} \le \zeta \le \zeta_{N}+1, ~N\ge 1.  
     \label{m1a456rcc} 
  \end{eqnarray}
  Substituting $\zeta=\zeta_{N}+1$ into (\ref{m1a456rcc}), and  using the property  $\lim_{\zeta \to \infty}(\zeta\rho)'= 0,$   we obtain
   \begin{eqnarray}
  (\zeta \rho)'{\Big |}_{\zeta=\zeta_{N}+1}  \ge   (\zeta \rho)'{\Big |}_{\zeta=\zeta_{N}}+ 
   {\frac {a \delta^{2}}{2M}} \ge {\frac {a \delta^{2}}{4M}}>0
  ~~{\rm when}~~N \gg 1, 
     \label{m1a456rccd} 
  \end{eqnarray}
  contradicting  the fact that  $\lim_{\zeta \to \infty}(\zeta\rho)'= 0.$ Thus, 
   property (\ref{m1rrfd2b}) cannot hold. It remains to assume, for contradiction,  that  the first possibility in  (\ref{m1rrfd2ac}) occurs,  i.e. 
           \begin{eqnarray}
  -M \le    \zeta_{N}\rho(\zeta_{N})\theta'(\zeta_{N}) \le -2\delta ~~\forall N \ge 1. 
   \label{m1rrfd2gga}
  \end{eqnarray}
The first step in obtaining a contradiction to (\ref{m1rrfd2gga}) is to assume that $\hat \zeta$ is chosen to satisfy (\ref{m1a34tt}) and 
(\ref{m1ad1tgggg14}), and also
         \begin{eqnarray}
   \hat \zeta >  {\frac {4M^{2}}{a\delta^{2}}}.
   \label{m1rrr}
  \end{eqnarray}
Next, we combine (\ref{m1art}) with the lower bound  ${\frac { \left ( \zeta \rho \right )' }{\rho}} \ge -{\frac {\delta^{2}}{aM^{2}}}$ 
  in (\ref{m1ad1tgggg14}),   
and the fact that  $ 0\le \zeta\rho\le M~~\forall \zeta \ge 0,$ and  obtain
   \begin{eqnarray}
     \left ( \zeta^{2}\rho^{2}\theta' \right )'\le    \delta^{2}~~\forall \zeta > \hat \zeta.
  \label{m1artabbc}
  \end{eqnarray}
  It follows from (\ref{m1a34}), (\ref{m1rrfd2gga}) and an integration of (\ref{m1artabbc}) that, for each $N \ge 1,$ 
     \begin{eqnarray}
     \zeta^{2}\rho^{2}\theta'  \le  \zeta_{N}^{2} \rho^{2} (\zeta_{N}) \theta'(\zeta_{N})  + \delta^{2}(\zeta-\zeta_{N})
    \le -2\delta^{2}+ \delta^{2} = -\delta^{2},~~\zeta_{N}\le \zeta \le  \zeta_{N}+1.
  \label{m1artabbc1}
  \end{eqnarray}
    Thus, since $0 < \zeta \rho \le M~~\forall \zeta >0,$ it follows from (\ref{m1artabbc1}) that 
  \begin{eqnarray}
     \zeta \rho \theta' \le   
     -{\frac {\delta^{2}}{\zeta \rho}}\le  -{\frac {\delta^{2}}{M}} ,~~\zeta_{N}\le \zeta \le  \zeta_{N}+1, ~N \ge 1. 
  \label{m1artabbc1a}
  \end{eqnarray}
 To make use of these properties, we write equation (\ref{m1aerddd1}) as 
 \begin{eqnarray}
  (\zeta \rho)'' =\zeta \rho \theta' \left (  \theta'  +a\zeta \right ) + \zeta \rho  -\zeta \rho^{3}.  
     \label{m1ar44rrrr}
  \end{eqnarray}
Combining (\ref{m1aeez}), (\ref{m1rrr}), (\ref{m1artabbc1})  and the fact that $0 \le \zeta \rho \le M~~\forall \zeta \ge 0,$ 
with equation (\ref{m1ar44rrrr}), we
obtain 
 \begin{eqnarray}
  (\zeta \rho)'' \le -{\frac {\delta^{2}}{M}} \left (  {\frac {a\zeta}{2}} \right ) + M \le
  -{\frac {\delta^{2}a\zeta}{4M}},  ~~\zeta_{N}\le \zeta\le \zeta_{N}+1, N \ge 1.    
     \label{m1ar44r32}
  \end{eqnarray}
Integrating (\ref{m1ar44r32}) from $\zeta_{N}$ to $\zeta$ gives
 \begin{eqnarray}
    (\zeta \rho)'  \le   (\zeta \rho)'{\Big |}_{\zeta=\zeta_{N}} -{\frac {\delta^{2}a}{4M}}
    \left ( {\frac {\zeta^{2}}{2}}-{\frac {\zeta_{N}^{2}}{2}} \right ),~~\zeta_{N}\le \zeta\le \zeta_{N}+1, N \ge 1.
     \label{m1r}
  \end{eqnarray}
Setting $\zeta=\zeta_{N}+1$ in (\ref{m1r}), and using the fact $(\zeta \rho)'{\Big |}_{\zeta=\zeta_{N}} \to 0$ as $N \to \infty,$ 
  we obtain
\begin{eqnarray}
    (\zeta \rho)'{\Big |}_{\zeta=\zeta_{N}+1}  \le   (\zeta \rho)'{\Big |}_{\zeta=\zeta_{N}} -{\frac {\delta^{2}a}{4M}}<
     -{\frac {\delta^{2}a}{8M}}<0
  ~~{\rm when}~~N \gg 1,
     \label{m1dr44}
  \end{eqnarray}
contradicting the fact that $(\zeta \rho)'{\Big |}_{\zeta=\zeta_{N}+1} \to 0$ as $N\to \infty.$ 
We conclude that   (\ref{m1rrfd2gga}) does not hold, hence  $E(\infty)=0,$ as claimed.  This completes the proof of Lemma~\ref{lem3a}.

 \smallskip \noindent Our final   result is 

       \begin{lem}  \label{lem7}   
      Let $\rho_{0}>0$ and $a>0.$  Then  
      the     solution of (\ref{m1aerddd1})-(\ref{m1aerddd2})-(\ref{m1aerddd3}) satisfies
                \begin{eqnarray}
 \lim_{\zeta \to \infty}\theta'(\zeta)=0~~{\rm and}~~~~\lim_{ \zeta \to \infty} \zeta\rho(\zeta) ={\sqrt {-E(\infty)}}.
  \label{m1a1}
  \end{eqnarray}   
    \end{lem}
    
    \smallskip \noindent {\bf Remarks.} {\bf (i)} Proving the first property in completes the
    proof of   (\ref{eq:d22}), which in turn completes the
    proof of    Theorem~\ref{th1aaa}.
  
  \noindent   
    {\bf (ii)} Much  of the proof    uses  the same basic approach  as in the proof of Lemma~\ref{lem3a}. 
     Due to the importance of  Lemma~\ref{lem7}   we    give complete details.

 \smallskip \noindent    
 {\bf Proof.} First,  recall  from   (\ref{m123jww1}) in Lemma~\ref{lem2} that $M>0$ exists such that
  \begin{eqnarray}
  \label{m123jww1a}
  0 \le \zeta \rho \le M~~\forall \zeta \ge 0.
    \end{eqnarray}
 We also  recall   from (\ref{m1aqq})-(\ref{m1adedtryaac}) that  $E$ 
 satisfies
  $ E(0)=\rho_{0}^{2}~~{\rm and}~~ E' <0~~\forall \zeta > 0.$
These  properties and  Lemma~\ref{lem3a} imply that $\zeta_{\rho_{0}}>0$ exists such that
  \begin{eqnarray}
  E'(\zeta)<0~\forall \zeta \in (0,\zeta_{\rho_{0}})~~{\rm and}~~E(\zeta_{\rho_{0}})=0,
 \label{m134}
  \end{eqnarray}
     \begin{eqnarray}
  E<0~~{\rm and}~~E'<0~~\forall \zeta > \zeta_{\rho_{0}}~~{\rm and}~~-\infty<E(\infty)<0.
  \label{m15}
  \end{eqnarray} 
  We conclude from (\ref{m15}) 
   that   $\lambda_{0}>0$ exists such that
 \begin{eqnarray}
  E(\zeta)= \left ( \left ( \zeta \rho  \right )' \right )^{2}+
 \left ( \zeta \rho  \right )^{2} \left ( (\theta')^{2}-1   +{\frac {\rho^{2}}{2}} \right )<-\lambda_{0}<0~~\forall
  \zeta >
 \zeta_{0}+1.  
  \label{m1ae}
  \end{eqnarray}
   It follows from (\ref{m1ae}) that  
   \begin{eqnarray}
   \left( \theta'(\zeta) \right)^{2}  <1~~  \forall \zeta >\zeta_{0}+1.
  \label{m1adde4}
  \end{eqnarray}
Also, we conclude from (\ref{m123jww1a}), (\ref{m1ae}) and (\ref{m1adde4}) that $ \liminf_{\zeta \to \infty} \zeta\rho(\zeta)>0.$ 
 From this  property and (\ref{m123jww1a}) it follows     that   $m \in (0,M]$ exists such  such that
           \begin{eqnarray}
   0<  m\le \zeta\rho\le M~~\forall \zeta > \zeta_{0}+1.
  \label{m1aedd3}
  \end{eqnarray}

 \noindent Our next goal is to  make use of  these properties to prove that 
      \begin{eqnarray}
   \lim_{\zeta \to \infty} \left ( \theta'(\zeta)\right )^{2}=0.
  \label{m1ad67tr1a}
  \end{eqnarray}
We assume, for contradiction, that  (\ref{m1ad67tr1a})  doesn't   hold, hence  
$ \limsup_{\zeta \to \infty} \left ( \theta'(\zeta)\right )^{2}>0.$
    This property
  and the fact that $0 \le \left ( \theta'(\zeta)\right )^{2}<1$ 
  imply that   $\delta \in (0,1) $ exists, and also a positive, increasing, unbounded sequence $ \left (\zeta_{N} \right )$ exists   such that
      \begin{eqnarray}
   \delta^{2} \le \left (  \theta' (\zeta_{N} )\right )^{2}<1~~\forall N \ge 1.
    \label{m1rrfd2c}
   \end{eqnarray} 
 Therefore, by considering subsequences if necessary,  either 
       \begin{eqnarray}
  -1  <    \theta'(\zeta_{N}) \le - \delta ~~\forall N \ge 1~~{\rm or}~~ \delta \le 
    \theta'(\zeta_{N}) <1~\forall N \ge 1.
  \label{m1rrfd2a}
  \end{eqnarray}
The proof of  (\ref{m1ad67tr1a}) is complete if we obtain a  
  contradiction to each     case  in (\ref{m1rrfd2a}). For this
  we again    use    three basic properties of solutions. 
First, because of (\ref{m123jww1a}) and the property $\lim_{\zeta \to \infty}\zeta_{N}=\infty,$ 
we can assume that
 $\zeta_{1}> \zeta_{0}+1$   is large enough so that 
      \begin{eqnarray}
  0<\rho(\zeta )<1~~\forall \zeta \ge \zeta_{1}.
  \label{m1ad67dddtr11}
  \end{eqnarray}
 It  follows from    (\ref{m1ad1tgggg}) in Lemma~\ref{lem6} 
 and the property $\lim_{\zeta \to \infty}\zeta_{N}=\infty,$ that 
 we can also assume that  
   $\zeta_{1}>\zeta_{0}+1 $  is large enough so that
  \begin{eqnarray}
  -{\frac {\delta}{2a}}\left (   {\frac {m}{M}}    \right )^{4}<
  {\frac {\left (\zeta \rho \right )'}{\rho}}<{\frac {\delta}{2a}}\left (   {\frac {m}{M}}    \right )^{4}~~\forall \zeta \ge \zeta_{1},
   \label{m1addder34rddeedd2}
  \end{eqnarray}
  \begin{eqnarray}
  -{\frac {a\delta m^{3}}{4M^{2}}} <  
  \left (\zeta \rho \right )' < {\frac {a\delta m^{3}}{4M^{2}}}~~\forall \zeta \ge \zeta_{1},
   \label{m1addder34rddeedd112}
  \end{eqnarray}
     Suppose, now, that  the second case in (\ref{m1rrfd2a})  occurs, i.e. 
            \begin{eqnarray}
0< \delta \le 
    \theta'(\zeta_{N}) <1~\forall N \ge 1.
  \label{m1rrfd2aqq}
  \end{eqnarray} 
 Again, we    make use of the equation 
    \begin{eqnarray}
     \left ( \zeta^{2}\rho^{2}\theta' \right )'= -a \left (\zeta \rho \right )^{2}
     {\frac { \left (   \zeta \rho   \right )'}{\rho}}.
  \label{m1art2}
  \end{eqnarray}
 Substituting the upper bounds $(\zeta \rho)^{2} \le M^{2}$ and 
$ {\frac {\left (\zeta \rho \right )'}{\rho}}<{\frac {\delta}{2a}}\left (   {\frac {m}{M}}  
  \right )^{4}$ into (\ref{m1art2}) gives
    \begin{eqnarray}
     \left ( \zeta^{2}\rho^{2}\theta' \right )'\ge - {\frac {\delta m^{4}}{2M^{2}}},~~\zeta \ge \zeta_{1}.
  \label{m58erddddd1}
 \end{eqnarray}
Integrating  (\ref{m58erddddd1}) from $\zeta_{N}$ to $\zeta,$ and using $\theta'(\zeta_{N})\ge \delta,$ we get
    \begin{eqnarray}
     \zeta^{2}\rho^{2}\theta'   \ge 
     \zeta_{N}^{2}\rho^{2}(\zeta_{N}) \delta
         -{\frac {\delta  m^{4}}{2M^{2}}}(\zeta-\zeta_{N}),~~\zeta \ge \zeta_{N},~~N\ge 1.
  \label{m586dd1}
 \end{eqnarray}
Dividing   (\ref{m586dd1}) by $\zeta^{2}\rho^{2},$ and   using   the fact that
 $m^{2} \le \zeta^{2}\rho^{2}  \le M^{2},$   we obtain  
 \begin{eqnarray}
    \theta'   \ge 
    {\frac {\zeta_{N}^{2}\rho^{2}(\zeta_{N})}{ \zeta^{2}\rho^{2}}} \delta -
         {\frac {\delta m^{4}}{2\zeta^{2}\rho^{2}M^{2}}}\ge  
	{\frac {\delta m^{2}}{2M^{2}}}, ~~\zeta_{N} \le \zeta \le \zeta_{N}+1,~~N\ge 1.
  \label{m586dddddd}
 \end{eqnarray}
 Next, we again make use of the equation  
    \begin{eqnarray}
  \left (\zeta \rho\right )'' =\zeta \rho \left ( (\theta')^{2}+a\zeta\theta'+1 -\rho^{2}\right ).
     \label{mderddd11}
 \end{eqnarray}
  Recall from(\ref{m1aedd3})  that $m \le \zeta \rho \le M$ when $\zeta \ge \zeta_{0}+1$, and that $\zeta_{N} \gg 1 $ when  $N \gg 1.$ 
 Combining these properties  with  (\ref{m1ad67dddtr11}),  (\ref{m586dddddd})  and equation  (\ref{mderddd11}),
 we obtain
    \begin{eqnarray}
  \left (\zeta \rho\right )'' \ge \zeta \rho \left (  a\zeta\theta' \right )
  \ge {\frac {a \delta m^{3}}{2M^{2}}},~~\zeta_{N} \le \zeta \le \zeta_{N}+1,~~N\gg 1.
     \label{mderdddd1q}
 \end{eqnarray}
   Integrating (\ref{mderdddd1q}) from $\zeta_{N}$ to $\zeta,$
 and making use of the lower bound in (\ref{m1addder34rddeedd112}), we get
    \begin{eqnarray}
  \left (\zeta \rho\right )'  
  \ge  -{\frac {a \delta m^{3}}{4M^{2}}}   + {\frac {a \delta m^{3}}{2M^{2}}}\left (\zeta -\zeta_{N}
  \right ),~~\zeta_{N} \le \zeta \le \zeta_{N}+1,~~N\gg 1.
     \label{mded765q}
 \end{eqnarray}
 Thus,
     \begin{eqnarray}
  \left (\zeta \rho\right )'{\Big |}_{\zeta=\zeta_{N}+1}  
  \ge   {\frac {a \delta m^{3}}{4M^{2}}}>0,~~\zeta_{N} \le \zeta \le \zeta_{N}+1,~~N\gg 1,
     \label{mded7d65qq}
 \end{eqnarray} 
  contradicting  the fact that $\lim_{\zeta \to \infty}\left (\zeta \rho\right )'=0.$
 Thus,       (\ref{m1rrfd2aqq}) cannot hold,   eliminating the  
 second case in (\ref{m1rrfd2a}). 
     Next, suppose   that  the first case in (\ref{m1rrfd2a})  occurs, i.e. 
            \begin{eqnarray}
  -1<  \theta'(\zeta_{N}) \le -\delta~\forall N \ge 1.
  \label{m1rrfd2aqw}
  \end{eqnarray}
  Substituting    $(\zeta \rho)^{2} \le M^{2}$ and 
$ {\frac {\left (\zeta \rho \right )'}{\rho}}>-{\frac {\delta}{2a}}\left (   {\frac {m}{M}}  
  \right )^{4}$ into  (\ref{m1art2})  gives
    \begin{eqnarray}
     \left ( \zeta^{2}\rho^{2}\theta' \right )'\le  {\frac {\delta m^{4}}{2M^{2}}},~~\zeta \ge \zeta_{1}.
  \label{m58ertre1}
 \end{eqnarray} 
 Integrating   (\ref{m58ertre1}) from $\zeta_{N}$ to $\zeta,$ and using $\theta'(\zeta_{N}) \ge \delta,$ we get
    \begin{eqnarray}
     \zeta^{2}\rho^{2}\theta'   \le
     -\zeta_{N}^{2}\rho^{2}(\zeta_{N}) \delta+
         {\frac {\delta m^{4}}{2M^{2}}}(\zeta-\zeta_{N}),~~\zeta_{N} \le \zeta \le \zeta_{N}+1,~~N\gg 1.
  \label{m431}
 \end{eqnarray}
 Divide  (\ref{m431}) by $\zeta^{2}\rho^{2},$   make  use  the  bound  $m \le \zeta \rho \le M$  
  and obtain
    \begin{eqnarray}
     \theta'   \le
     - {\frac {m^{2}\delta}{M^{2}}}+
         {\frac {\delta m^{2}}{2M^{2}}}(\zeta-\zeta_{N})\le- {\frac {m^{2}\delta}{2M^{2}}},~~\zeta_{N} \le \zeta \le \zeta_{N}+1,~~N\ge 1.
  \label{m4rrr1}
 \end{eqnarray}
 Next, recall from  (\ref{m1adde4})   that $(\theta')^{2} < 1~~\forall \zeta > \zeta_{0}.$
 From this property, (\ref{m4rrr1}) and the fact that $\zeta_{N} \to \infty$ as $ N \to \infty,$ we conclude that
     \begin{eqnarray}
     (\theta')^{2}+a\zeta\theta'+1   \le
     - {\frac {3am^{2}\delta}{8M^{2}}},~~\zeta_{N} \le \zeta \le \zeta_{N}+1,~~N\gg 1.
  \label{m4drrr11}
 \end{eqnarray}
 Combining (\ref{mderddd11}) 
  with (\ref{m4drrr11}) and the lower bound $\zeta\rho \ge m$ when $\zeta \ge \zeta_{0}+1$
  gives
   \begin{eqnarray}
  \left (\zeta \rho\right )'' \le  -
   {\frac { 3a \delta m^{3}}{8M^{2}}},~~\zeta_{N} \le \zeta \le \zeta_{N}+1,~~N\gg 1.
     \label{mdedddrdddd11}
 \end{eqnarray}
 It follows from an integration of (\ref{mdedddrdddd11}) and the upper bound in (\ref{m1addder34rddeedd112}) that
    \begin{eqnarray}
  \left (\zeta \rho\right )' \le  {\frac {a \delta m^{3}}{4M^{2}}}-
   {\frac {3 a \delta m^{3}}{8M^{2}}}(\zeta -\zeta_{N}),~~\zeta_{N} \le \zeta \le \zeta_{N}+1,~~N\gg 1.
     \label{mded989dddd12}
 \end{eqnarray}
 Finally, we conclude from (\ref{mded989dddd12}) that
  \begin{eqnarray}
  \left (\zeta \rho\right )'{\Big |}_{\zeta=\zeta_{N}+1}  
  \le  - {\frac {a \delta m^{3}}{8M^{2}}}<0,~~N\gg 1,
     \label{mddded7d653}
 \end{eqnarray} 
  contradicting the property $\lim_{N\to \infty}\left (\zeta \rho\right )'{\Big |}_{\zeta=\zeta_{N}+1}   =0.$
 This completes the proof that   $\lim_{\zeta \to \infty} \left (\theta'(\zeta)\right )^{2}=0,$
 Thus,  $\lim_{\zeta \to \infty} \theta'(\zeta) =0,$ and the the first property in (\ref{m1a1}) is proved.
    It remains to prove the second property in (\ref{m1a1}), namely
   \begin{eqnarray}
  \lim_{ \zeta \to \infty} \zeta\rho(\zeta) ={\sqrt {-E(\infty)}}.
 \label{m1add345}
  \end{eqnarray} 
 First, we write (\ref{m1aqq})
   as
    \begin{eqnarray}
  \left ( \zeta \rho  \right )^{2}
   \left ( (\theta')^{2}-1 \right ) =E-{\frac {1}{2}}\zeta^{2}\rho^{4}-\left ( \left ( \zeta \rho  \right )' \right )^{2}.
  \label{m1adedtryaaj}
  \end{eqnarray}
 Since $\left (\theta'\right )^{2}<1~~{\rm and}~~E<0~~{\rm when}~~\zeta>\zeta_{0},$ we can divide (\ref{m1adedtryaaj}) by  
 $(\theta')^{2}-1$  and get
  \begin{eqnarray}
  \left ( \zeta \rho  \right )^{2}
     ={\frac {E-{\frac {1}{2}}\zeta^{2}\rho^{4}-\left ( \left ( \zeta \rho  \right )' \right
   )^{2}}{(\theta')^{2}-1}}.
  \label{m1adedtryaajb}
  \end{eqnarray}
We conclude  from  the third property in (\ref{m1ad1tgggg}),  (\ref{mderddd11})  and (\ref{m1adedtryaajb}) 
  that
    \begin{eqnarray}
  \lim_{ \zeta \to \infty} \left (\zeta\rho(\zeta)\right )^{2} = -E(\infty).
 \label{m1add345q}
  \end{eqnarray}
 Property (\ref{m1add345}) follows from (\ref{m1add345q})  and the fact that $\zeta \rho >0~~\forall \zeta>0.$

\bigskip \noindent {\bf Competing interests and funding.}  
 There are no competing interests.


\begin{thebibliography}{99}
 
                  \bibitem{adk} {G. V. Arkrivis, V. A. Dougalis, O. A. Karakashian, W. R. McKinney}  {\em 
 Galerkin-finite element methods for the nonlinear Schrodinger equation,}
  Hellenice research in Mathematics and Informatics '92, Hellenic Math. Soc., (Athens, 1992),   421-442
 
      
      
      
      
                  \bibitem{budd} {C. J. Budd,  S. Chen and R. D. Russel,}  {\em 
 New self-similar solutions of the nonlinear Schrodinger equation with moving mesh computations,}
     J. Comp. Phys.
      {\bf 152} (1999), 756-789  
      
                      \bibitem{chai} {R. Y. Chiao, E. Garmire    and C. H. Townes}  {\em 
 Self-trapping of optical beams,}
  Phys. Rev. Lett.  {\bf 13} 
        (1964), 479-382  
   
   
                   \bibitem{chu} {P. K. Lu and L. XinPei,}  {\em 
 Low temperature plasma technology: methods and applications,}
   CRC Press
       (2013)   
   
   
   
   
   
   
    
     
    
    
    
    
    
    
    
                \bibitem{fib96} {G. Fibich}  {\em 
  Adiabatic law for self-focusing of optical beams,}
  Optics Letters, {\bf 21}    
        (1996), 1735-1737  
    
 
    
    
    
    
                   \bibitem{fp98} {G. Fibich and G. C. Papanicolau}  {\em 
 Self-focusing in the perturbed and unperturbed nonlinear Schrodinger equation in critical dimension,}
  SIAM J. Appl. Math 60 (1), 
        (1999), 183-240    
    
    
    
 
    
    
    
    
    
    
    
    
    
    
              \bibitem{ha} { A. Hasegawa,}   
{\em Optical solitons in fibers,} 
Springer-Verlag, Berlin  (1989) 
      
      
      
 
    
           \bibitem{has} { S. P. Hastings and J. B. McLeod,}   
{\em Classical Methods in Ordinary Differential Equations With Applications to Boundary Value Problems,} 
Amer. Math. Soc.  (2012)
  
  
     
    
    
            \bibitem{kop} {N. Kopell and M. Landman,}  {\em 
Spatial structure of the focusing  singularity of the nonlinear Schrodinger equation. A geometrical analysis.}
     SIAM J. Appl. Math.
      {\bf 55} (1995), 1297-1323  
    
    
           \bibitem{Land} {M. J. Landman,   G. C. Papanicolau, C. Sulem  and   P. Sulem,}  {\em 
 Rate of blowup for solutions of the nonlinear Schrodinger equation at critical dimension,}
    Phys. Rev. A. 
      {\bf 38} (1988), 3837-3843    
    
    
               \bibitem{langm} {I. Langmuir,}  {\em 
 Oscillations in ionized gases,}
    Proc. Nat. Acad. Sci. 
      {\bf 14 } (1928), 627-637 
    
    
                                     
              \bibitem{L1} {B. LeMesurier, G. C. Papanicolau, C. Sulem  and   P. Sulem,}  {\em 
 Focusing and multi-focusing solutions of the nonlinear Schrodinger equation,}
    Physica D 
      {\bf 31} (1988), 78-102 
    
    
      
      
           \bibitem{L} {B. LeMesurier, G. C. Papanicolau, C. Sulem  and   P. Sulem,}  {\em 
 Local structure of the self-focusing nonlinearity of the nonlinear Schrodinger equation,}
    Physica D 
      {\bf 32} (1988), 210-220        
      
     
              \bibitem{MCL} {D. W. McGlaughlin, G. C. Papanicolau, C. Sulem  and   P. Sulem,}  {\em 
 Focusing singularity of the cubic  nonlinear Schrodinger equation,}
    Phys. Rev. A. 
      {\bf 34} (1986), 1200-1210   
     
     
               \bibitem{new} {A. C. Newell,}  {\em 
 Solitons in Mathematics and Psysics,} in CBMS Applied Mathematical Series,
   SIAM 
          Phil. Pa.{\bf 48} (1985)  
     
     
      
     
     
     
     
                   \bibitem{ras} {J. J. Rasmussen and K. Rypdal}  {\em 
Blowup in nonlinear Schrodinger equations. 1. A general review and 2. Similarity structure of the blowup singularity,}
      Physica Scripta 
      {\bf 33} (1986),  481-504  
     
     
     
     
     
     
     
     
      
      
      
                \bibitem{rott} {V. Rottschafer and T. J. Kaper}  {\em 
Blowup in the nonlinear Schrodinger equation near critical dimension,}
     J. Math. Anal. Appl.
      {\bf 268} (2002),  517-549 
      
      
      
                  \bibitem{rott2} {V. Rottschafer and T. J. Kaper}  {\em 
Geometric theory for multi-bump self similar  blowup solutions of the cubic   nonlinear Schrodinger equation,}
     Nonlinearity
      {\bf 16} (3) (2003),  929-962    
  
      
        
      
       
         \bibitem{sul} { C. Sulem and P. Sulem,}   
{\em The Nonlinear Schrodinger Equation,} Springer , New York (1999)
   
 
 
 
 
 
 
 
                \bibitem{ts92} {Y. Tourigny and J. M. Sanz-Serna}  {\em 
The numerical study of blowup with application to a nonlinear Schrodinger equation,}
     J. Comp. Phys., 102
       (1992), 407-416 
 
 
 
 
 
 
 
 
 
 
 
 
 
 
 
 
 
 
 
 
 
    
    
         
                \bibitem{wang} {X. P. Wang} 
		 {\em 
On singular solutions of the nonlinear Schrodinger and Zahkarov equations,}
   Ph.D. Thesis, New York University, 1990
       
    
    
    
    
    
    
    
            \bibitem{zh1972} { V. E. Zakharhov,} {bf NEW}
	     		 {\em 
Collapse  of Langmuir waves,}    
Sov. Phys. JETP,  Vol. 55, no. 5 (1972), 908-914     
    
    
    
    
    
    
      
            \bibitem{z} { V. E. Zakharhov,} 
	     		 {\em 
Collapse and Self-focusing of Langmuir waves,}    
{\em Handbook of Plasma Physics,} Vol. 2, M. N. Rosenbluth and R. Z. Sagdeev, eds.,  Elsevier, Amsterdam (1984)   
  

          
 
 

\end{thebibliography}
 \end{document}